\documentclass[11pt, reqno]{amsart}
\usepackage{amsmath, amsthm, amscd, amsfonts, amssymb, graphicx, hyperref,color}
\newtheorem{theorem}{Theorem}[section]
\newtheorem{lemma}[theorem]{Lemma}

\newtheorem{corollary}[theorem]{Corollary}
\theoremstyle{definition}
\newtheorem{definition}[theorem]{Definition}
\newtheorem{example}[theorem]{Example}

\theoremstyle{remark}
\newtheorem{remark}[theorem]{Remark}
\numberwithin{equation}{section}
\begin{document}
\setcounter{page}{1}
\title[Right and left quotient of two bounded operators]{Right and left quotient of two bounded operators on Hilbert spaces}
\author{ Mohammed Benharrat}
\address{ Department of Mathematics and informatics, National Polytechnic School of Oran-Maurice Audin, BP 1523 Oran-El M'naouar, 31000, Oran, Algeria.}
\email{mohammed.benharrat@gmail.com;  benharrat@math.univ-lyon1.fr.}


\dedicatory{Dedicated to the memory of R. G. Douglas  (1938-2018)}
\date{09/11/2018.\\
	This work was supported by the Algerian research project: PRFU, no. C00L03ES310120180002.}
\begin{abstract}
We define a left  quotient as well as a right quotient of two bounded operators between Hilbert spaces, and we parametrize these two concepts using the Moore-Penrose inverse.  In particular, we show that the adjoint of a left quotient is a right quotient and conversely. An explicit formulae for computing  left (resp.   right)  quotient which correspond to adjoint, sum, and  product of given left (resp.   right) quotient of two bounded operators are also shown.
\end{abstract}

\maketitle
\section{Introduction }
Let $A$ and $B$ be bounded linear  operators on a Hilbert space $\mathcal{H}$ with \textbf{the kernel condition}
\begin{equation*}
\mathcal{N}(B)\subset \mathcal{N}(A).
\end{equation*}
 The quotient operator $[A/B]$ of $A$ and $B$ (not necessary bounded) is defined  as the mapping  $Bx\longmapsto Ax$, $x \in \mathcal{H}$. If
we write $\mathcal{G}(B,A)$ for the set $\{(Bx,Ax): x \in \mathcal{H}\}$ in the product Hilbert space
$\mathcal{H}\times \mathcal{H}$, then $\mathcal{G}(B ,A)$ is a graph and we can define $[A/B]$ as the operator
corresponding to this graph,  see \cite{Kaufman2}. A quotient (of bounded operators) so defined appeared  for the first time in the work of Dixmier \cite{Diximier49}  by the name  "op\'erateur $J$ uniforme" and investigated by  Kaufman under the name  "semi closed operators"  in \cite{Kaufman1, Kaufman2, Kaufman3}, were several characterizations are given. It is worth noting that this extension was anticipated well before by Douglas in  \cite{Douglas66}, where he used in the proof of his famous lemma. In \cite{Kaufman1} Kaufman  showed that a linear operator $T$ on $ \mathcal{H}$
is closed if and only if $T$ is represented as a quotient $[A/B]$ using $A$ and $B$ such that
$\mathcal{R}(A^{*})+\mathcal{R}(B^{*})=\left\{A^{*}x+B^{*}y : x; y \in \mathcal{H}\right\}$ is closed in $\mathcal{H}$, so that every closed operator
is included in the class of quotients ($\mathcal{R}(A)$ denotes the range of an operator $A$). Since many authors are studding  this class of operators. We can mention, for example, the study of algebraic operations,  adjoint, weak adjoint and selfadjoint extension of quotient operators are included in \cite{izumino89, izumino93, Hirasawa95,Hirasawa04,izuminoH01},  various  topologies on the sets of unbounded linear Hilbert space operators by means of quotient operators are investigated in \cite{Hirasawa07, Hirasawa11, Benharrat16,Kaufman4, Kittaneh}, and for other topics we can see \cite{izumino89HJ, Benharrat6, Benharrat14}. Recently,  Koliha in  \cite{Koliha}  extend Kaufman's results to operators between two Hilbert spaces $\mathcal{H}_1$  and $\mathcal{H}_2$  as follows: if $A$ is a bounded operator from $\mathcal{H}_1$ to $\mathcal{H}_2$ and $B$ is bounded operator on $\mathcal{H}_1$, we can define the quotient $[A/B]$ as the set $\mathcal{G}(B,A) =\{(Bx,Ax): x \in \mathcal{H}_1\}\subset \mathcal{H}_1 \times \mathcal{H}_2$. He showed that   $[A/B]$ is a linear operator from $\mathcal{H}_1$ to $\mathcal{H}_2$, with domain $\mathcal{R}(B)$ and range $\mathcal{R}(A)$ if and only if $\mathcal{N}(B)\subset \mathcal{N}(A)$. For more details see \cite{Koliha}. 

 This paper is a continuation and refinement of the research treatment of the class of quotient operators.  From now,  the concept of the  left quotient of bounded operators between two Hilbert spaces (shortly, a left quotient) means  the Douglas solution of a given range  operators inclusion. A right quotient means  quotients of bounded operators as defined above. In this work we extend this two definitions to bounded operators acting  between two possibly different Hilbert spaces. We also collect; unify and generalize earlier work and derive new or strengthened results, including a complete description of the this two class of operators   in terms of the Moore-Penrose pseudoinverse. We next prove that the adjoint of a left quotient is always exist and coincide with the right quotient of the adjoints of this given operators. We present reasonable left quotient which coincide with the adjoint of a given right quotient. We investigate a conditions under which the sum and product of two  left  (resp. right) quotient operators are also left  (resp. right) quotient operators and we show how to compute the sum and product of given quotients.  We also give some other applications which show the advantage of fractional representations of operators


\section{The left quotient of bounded operators}
In what follows $\mathcal{H}$, $\mathcal{H}_1$, $\mathcal{H}_2$, etc.  denote Hilbert spaces endowed with the appropriate  scalar product  and the associated norm.  The inner product in $\mathcal{H}_1 \times \mathcal{H}_2$ is defined by $<(x,y),(x',y')>= <x,x'>+<y,y'>$.  For $T$  linear operator from $\mathcal{H}_1$ to $\mathcal{H}_2$, the symbols $\mathcal{D}(T)\subset \mathcal{H}_1$, $\mathcal{N}(T)\subset \mathcal{H}_1$ and $\mathcal{R}(T)\subset \mathcal{H}_2$ will denote the domain, null space and the range space of $T$, respectively. The set $\mathcal{G}(T)=\{(x,Tx): x\in\mathcal{D}(T)\} \subset \mathcal{H}_1 \times \mathcal{H}_2$ is called the graph of $T$. The operator $T$ is closed if and only if $\mathcal{G}(T)$ is a closed subset of  $\mathcal{H}_1 \times \mathcal{H}_2$, and is  densely defined if $ \overline{\mathcal{D}(T)}= \mathcal{H}_1$, where $\overline{\mathcal{D}(T)}$ denote the closure of $\mathcal{D}(T)$ in $\mathcal{H}_1$.  Denote by $\mathcal{B}(\mathcal{H}_1 , \mathcal{H}_2)$ the Banach space of all bounded linear operators from  $\mathcal{H}_1$ to $\mathcal{H}_2$. If $\mathcal{H}_1=\mathcal{H}_2$, write  $\mathcal{B}(\mathcal{H}_1 , \mathcal{H}_2)=\mathcal{B}(\mathcal{H}_1)$.   

Before to introduce the left quotient, we need the following theorem, which is referred to as the Douglas lemma \cite{Douglas66}. This lemma is a basic tool for our discussions.
\begin{theorem}\cite[Theorem 1.]{Douglas66} For  $A\in \mathcal{B}(\mathcal{H}_1, \mathcal{H}_3 )$ and $B\in \mathcal{B}(\mathcal{H}_2, \mathcal{H}_3 )$, the following statements
are equivalent:
\begin{enumerate}
	\item (\textbf{range inclusion}) $\mathcal{R}(A) \subseteq \mathcal{R}(B)$;
	\item (\textbf{majorization}) $AA^* \leq \lambda^2 BB^*$ for some $\lambda \geq 0$; and
	\item (\textbf{factorization}) there exists a bounded operator $C$ from $\mathcal{H}_1$ to $\mathcal{H}_2$ such that  $A = BC$.
\end{enumerate}
Moreover, if one of this conditions holds, then there exists a unique operator $C$ so that
\begin{itemize}
	\item [(a)] $\left\|C\right\|^2=\inf\{\mu : AA^* \leq \mu BB^*\}$;
	\item [(b)] $ \mathcal{N}(C)=\mathcal{N}(A)$; and
	\item [(c)] $\mathcal{R}(C) \subseteq \overline{\mathcal{R}(B^*)}$
\end{itemize}
We shall call this uniquely determined $C$ the Douglas solution of the equation $BX=A$.
\end{theorem}
Let us point out that the operator $C$ is uniquely determined by the the condition (c) and verifies the two conditions (a) and (b); in general the  operator $C$ may be not unique but (a) and (b) are valid, as shown in the following examples:
\begin{example}\label{examp1}
Let 
$$ A=\begin{bmatrix} 0 & 1\\ 0 & 0\end{bmatrix}, \qquad \text{ and } \qquad  B=\begin{bmatrix} 1 & 1\\ 0 & 0\end{bmatrix}.$$
We have $\mathcal{R}(A) = \mathcal{R}(B)=span\{e_1\}$, $AA^* = \frac{1}{2} BB^*$ and  $A = CB$ for 
$$ C=\begin{bmatrix} \alpha & \beta \\ -\alpha & 1-\beta \end{bmatrix}\qquad \text{ for  }  \alpha , \beta \in \mathbb{R}.  $$
Now if we assume that the matrix $C$ such that $ \mathcal{N}(C)=\mathcal{N}(A)$ (the condition (b)), we get $\alpha =0$ and $A = BC$ is verified for all $C$, with
$$ C=\begin{bmatrix} 0 & \beta \\ 0 & 1-\beta \end{bmatrix}\qquad \text{ for all }  \beta \in \mathbb{R}.  $$
Further, if we added the condition (c), we obtain $\beta=\frac{1}{2}$; so
$$ C=\begin{bmatrix} 0 & \frac{1}{2} \\ 0 & \frac{1}{2} \end{bmatrix}.  $$
\end{example}
\begin{example}\label{examp2}
Let 
$$ A=\begin{bmatrix} 0 & 0& 1\\ 0 & 0 & 0\\ 0 & 0& 0\end{bmatrix}, \qquad \text{ and } \qquad  B=\begin{bmatrix} 0 & 0& 1\\ 0 & 1&0\\ 0&0&0\end{bmatrix}.$$
By a  simple calculus, we get  $\mathcal{R}(A)=span\{e_1\} \subset \mathcal{R}(B)=span\{e_1, e_2\}$, $AA^* \leq  BB^*$ and  $A = BC$ for 
$$ C=\begin{bmatrix} \alpha & \beta &\gamma \\ 0 & 0&0\\0&0&1 \end{bmatrix}\qquad \text{ for  all }  \alpha , \beta,\gamma \in \mathbb{R}.  $$
Now if we assume that the matrix $C$ such that $ \mathcal{N}(C)=\mathcal{N}(A)=span\{e_1, e_2\}$ (the condition (b)), we get $\alpha =\beta=0$ and $A = BC$ is verified for all $C$, with
$$ C=\begin{bmatrix} 0 & 0 &\gamma\\ 0 & 0&0\\0&0&1 \end{bmatrix}\qquad \text{ for all }  \gamma\in \mathbb{R}.  $$
Furthermore; $\left\|C\right\|=1$; at this stage $C$ is not  unique. Hence $C$ is uniquely determined by using the condition (c), witch gives
$$ C=\begin{bmatrix} 0 & 0 &0\\ 0 & 0&0\\0&0&1 \end{bmatrix}.$$
\end{example}
In view of the observation above, a definition of a left quotient of two operators can be given as follows:
\begin{definition}[\textbf{Left quotient of bounded operators}]Let  $A\in \mathcal{B}(\mathcal{H}_1, \mathcal{H}_3 )$ and $B\in \mathcal{B}(\mathcal{H}_2, \mathcal{H}_3 )$ such that  the \textbf{range inclusion} : $\mathcal{R}(A) \subseteq \mathcal{R}(B)$ holds. The unique operator 
$$C: \mathcal{H}_1 \longrightarrow \mathcal{N}(B)^{\bot}$$
as a Douglas solution of $A=BX$ is called \textbf{left quotient } of $A$ by $B$, and will be noted by $[B \backslash A]$.
\end{definition}
A direct consequence of  the Douglas lemma, we have:
\begin{itemize}
	\item [(i)]  $[B \backslash A]$ is bounded operator from  $\mathcal{H}_1$ to $\mathcal{H}_2$, with $\left\|[B \backslash A]\right\|=\sqrt{\inf\{\mu : AA^* \leq \mu BB^*\}}$. 
	\item [(ii)] $ \mathcal{N}([B \backslash A])=\mathcal{N}(A)$ and $[B \backslash A]$ from $\mathcal{N}(A)^{\bot}$ to $\mathcal{N}(B)^{\bot}$ is an injective operator.
	\item [(iii)] $B[B \backslash A]x=Ax$ for all $x\in \mathcal{H}_1$. So, the left quotient $[B \backslash A]$ is the unique solution of the equation $A=BX$.
\end{itemize}
Next, in order to  give a complete description of the left  quotient using the Moore-Penrose inverse, let recall  the definition of this generalized inverse.
\begin{definition}\cite{israel}
Let $T\in \mathcal{B}(\mathcal{H}_1, \mathcal{H}_2 )$. Then there exists a unique closed densely defined
operator $T^{\dag}$,   with domain $\mathcal{D}(T^{\dag})= \mathcal{R}(T) \oplus \mathcal{R}(T)^\bot$ and has the
following properties:
\begin{enumerate}
	\item $T^{\dag}Tx = P_{\mathcal{N}(T)^{\bot}} x$, for all $x \in \mathcal{H}_1$.
	\item $TT^{\dag}y = P_{\overline{\mathcal{R}(T)}}y$, for all $y \in \mathcal{D}(T^{\dag})$.
	\item $\mathcal{N}(T^{\dag})=  \mathcal{R}(T)^\bot$.
\end{enumerate}
Where $P_{\mathcal{M}}$ denotes the orthogonal projection onto a closed subspace $\mathcal{M}$.
\end{definition}
This unique operator $T^{\dag}$ is called the Moore-Penrose inverse of $T$.  (or  the Maximal
Tseng generalized Inverse in the terminology of \cite{israel}). In this case, $T$ is written in a matrix form with respect to mutually orthogonal subspaces decompositions as follows
$$T =\begin{bmatrix}
T_1& 0 \\
 0& 0
\end{bmatrix} : \begin{bmatrix}
\mathcal{N}(T)^{\bot}\\
 \mathcal{N}(T) 
\end{bmatrix}\longrightarrow  \begin{bmatrix}
 \overline{\mathcal{R}(T)}\\
\mathcal{R}(T)^{\bot}
\end{bmatrix} 
 $$
with $T_1$ is an invertible operator from $\mathcal{N}(T)^{\bot}$ to $\mathcal{R}(T)$. Hence, its Moore-Penrose inverse is given by
$$  T^{\dag} =\begin{bmatrix}
 T_1^{-1}& 0\\
 0 & 0
\end{bmatrix}: \begin{bmatrix}
 \mathcal{R}(T)\\
\mathcal{R}(T)^{\bot}
\end{bmatrix}\longrightarrow  \begin{bmatrix}
 \mathcal{N}(T)^{\bot}\\
 \mathcal{N}(T)
\end{bmatrix} .$$
It is well known  that $T^{\dag}$ is bounded if and only if  $\mathcal{R}(T)$ is closed in $\mathcal{H}_2 $.
It can seen that if $T$ is an injective then $T^{\dag}y=T^{-1}y$ for all $y\in \mathcal{R}(T)$.

In the next result, we redefine  the left quotient via the Moore-Penrose inverse. 
\begin{lemma}\label{lqd} Let  $A\in \mathcal{B}(\mathcal{H}_1, \mathcal{H}_3 )$ and $B\in \mathcal{B}(\mathcal{H}_2, \mathcal{H}_3 )$ such that $\mathcal{R}(A) \subseteq \mathcal{R}(B)$. Then $$[B \backslash A]=B^{\dag}A:  \mathcal{H}_1 \longrightarrow  \mathcal{N}(B)^{\bot}.$$ In particular, if $B$ is invertible then $[B \backslash A]=B^{-1}A$.
\end{lemma}
\begin{proof} Since  $\mathcal{R}(A) \subseteq \mathcal{R}(B)$,  $B^{\dag}A$ is bounded operator from $\mathcal{H}_1$ to $\mathcal{H}_2$ as a composition of closed operator and bounded one (see \cite[Problem 5.22, p 167]{Kato} ). Further,  $A$ and $B$  has the following form
$$A =\begin{bmatrix}
A_1& 0 \\
 0& 0
\end{bmatrix} : \begin{bmatrix}
\mathcal{N}(A)^{\bot}\\
 \mathcal{N}(A) 
\end{bmatrix}\longrightarrow  \begin{bmatrix}
 \mathcal{R}(B)\\
\mathcal{R}(B)^{\bot}
\end{bmatrix} 
 $$
and
$$  B =\begin{bmatrix}
 B_1& 0\\
 0 & 0
\end{bmatrix}: \begin{bmatrix}
 \mathcal{N}(B)^{\bot}\\
 \mathcal{N}(B)
\end{bmatrix}\longrightarrow  \begin{bmatrix}\mathcal{R}(B)\\
\mathcal{R}(B)^{\bot}
\end{bmatrix},$$
with $B_1$ is an invertible operator from $\mathcal{N}(B)^{\bot}$ to $\mathcal{R}(B)$. Moreover,
$$  B^{\dag} =\begin{bmatrix}
 B_1^{-1}& 0\\
 0 & 0
\end{bmatrix}: \begin{bmatrix}
 \mathcal{R}(B)\\
 \mathcal{R}(B)^{\bot}
\end{bmatrix}\longrightarrow  \begin{bmatrix}
 \mathcal{N}(B)^{\bot}\\
 \mathcal{N}(B)
\end{bmatrix} .$$
 Then 
$$B^{\dag}A= \begin{bmatrix}
 B_1^{-1}A_1& 0\\
 0 & 0
\end{bmatrix}: \begin{bmatrix}
\mathcal{N}(A)^{\bot}\\
 \mathcal{N}(A) 
\end{bmatrix}\longrightarrow  \begin{bmatrix}
 \mathcal{N}(B)^{\bot}\\
 \mathcal{N}(B)
\end{bmatrix} .$$
This implies that
\begin{equation}\label{mpl1}
BB^{\dag}A= \begin{bmatrix}
 A_1& 0\\
 0 & 0
\end{bmatrix}=A: \begin{bmatrix}
 \mathcal{N}(B)^{\bot}\\
 \mathcal{N}(B)
\end{bmatrix}\longrightarrow   \begin{bmatrix}
 \mathcal{R}(B)\\
\mathcal{R}(B)^{\bot}
\end{bmatrix}.\end{equation}

On the other hand, 
\begin{equation}\label{mpl2}
BB^{\dag}A=P_{\overline{\mathcal{R}(B)}}Ax=Ax \qquad \text{ for all } x\in \mathcal{H}_1.
\end{equation}
By \eqref{mpl1} and  \eqref{mpl2}, it follows that  $BB^{\dag}A=A$ with
$$B^{\dag}A: \mathcal{H}_1 \longrightarrow  \mathcal{N}(B)^{\bot}.$$ Now by the uniqueness of the Douglas solution we have the desired result.
\end{proof}
\begin{remark}For $A\in \mathcal{B}(\mathcal{H})$, we have
	\begin{itemize}
		\item $[I \backslash A]=A$.
		\item If $A$ is surjective, then $[A \backslash I]$  is the right inverse of $A$ and given by $[A \backslash I]=A^{\dag}=A^*(AA^*)^{-1}$.
	\end{itemize}
\end{remark}
\begin{corollary}Let $A\in  \mathbb{C}^{n\times m}$  and $B\in \mathbb{C}^{n\times p}$ such that $\mathcal{R}(A) \subseteq \mathcal{R}(B)$. Then $[B \backslash A]=B^{\dag}A\in \mathbb{C}^{m\times p}$.
\end{corollary}
\begin{example} For the matrices $A$ and $B$ given in Example \ref{examp1}  and Example  \ref{examp2}, the left quotient of $A$ by $B$ is given by

$$ [B \backslash A]=B^{\dag}A=\begin{bmatrix} 0 & \frac{1}{2} \\ 0 & \frac{1}{2} \end{bmatrix} $$
and
$$ [B \backslash A]=B^{\dag}A=\begin{bmatrix} 0 & 0 &0\\ 0 & 0&0\\0&0&1 \end{bmatrix}$$
respectively. 
\end{example}
If $A\in \mathcal{B}(\mathcal{H}_1, \mathcal{H}_3 )$ and $B\in \mathcal{B}(\mathcal{H}_2, \mathcal{H}_3 )$ such that $\mathcal{R}(A) = \mathcal{R}(B)$, then we can define $[B \backslash A]$ and $[A \backslash B]$ as follows
$$[B \backslash A]: \mathcal{H}_1 \longrightarrow \mathcal{N}(B)^{\bot}$$
and
$$[A \backslash B]: \mathcal{H}_2 \longrightarrow \mathcal{N}(A)^{\bot}.$$
with $$[B \backslash A][A \backslash B] \text{  is the identity operator on }
\mathcal{N}(B)^{\bot} $$ and 
$$[A \backslash B][B \backslash A] \text{  is the identity operator on }
\mathcal{N}(A)^{\bot}.$$
 Therefore, $[B \backslash A]$ is an invertible operator from  $\mathcal{N}(A)^{\bot}$ to $\mathcal{N}(B)^{\bot}$ with
$$ [B \backslash A]^{-1} =[A \backslash B].$$
\begin{corollary}\label{corl1}Let $A\in \mathcal{B}(\mathcal{H}_1, \mathcal{H}_3 )$ and $B\in \mathcal{B}(\mathcal{H}_2, \mathcal{H}_3 )$  have the same range. Then $[B \backslash A]$ is an invertible operator from  $\mathcal{N}(A)^{\bot}$ to $\mathcal{N}(B)^{\bot}$ with
$$ [B \backslash A]^{-1} =[A \backslash B] : \mathcal{N}(B)^{\bot} \longrightarrow \mathcal{N}(A)^{\bot}$$
and bounded.
\end{corollary}

\section{The right quotient of bounded operators}
We now present the definition of the second kind of quotient of two bounded operators in the following generalized form.
\begin{definition}[\textbf{Right quotient of bounded operators}]Let  $A \in \mathcal{B}(\mathcal{H}_1, \mathcal{H}_2 )$ and  $ B\in \mathcal{B}(\mathcal{H}_1, \mathcal{H}_3 )$ such that the \textbf{kernel inclusion}: $\mathcal{N}(B) \subseteq \mathcal{N}(A)$ holds.  \textbf{The right quotient operator} $[A/B]$ of $A$ by $B$  is defined  as the mapping  $$Bx\longmapsto Ax; \quad x \in \mathcal{H}_1.$$ 
\end{definition}
A direct consequence of  this definition, we have:
\begin{itemize}
	\item [(i)] If we write $\mathcal{G}(B,A)$ for the set $\{(Bx,Ax): x \in H_1\}$ in the product Hilbert space
$\mathcal{H}_3\times \mathcal{H}_2$, then $\mathcal{G}(B ,A)$ is a graph and we can define $[A/B]$ as the operator from $\mathcal{H}_3$ to $\mathcal{H}_2$
corresponding to this graph.
\item [(ii)] $[A/ B]Bx=Ax$ for all $x\in \mathcal{H}_1$. So, the right quotient $[A/B]$ is the unique solution of the equation $A=XB$.
\item [(iii)] $ \mathcal{D}([A /B])=\mathcal{R}(B)$,  $\mathcal{R}([A/ B])=\mathcal{R}(A)$ and $\mathcal{N}([A/ B])=B(\mathcal{N}(A))$.
\item [(iv)] If $\mathcal{R}(B)$ is closed, then $[B/A]$ is a bounded  operator from $\mathcal{R}(B)$ to $\mathcal{R}(A)$.
\item [(v)] According to the decomposition $\mathcal{H}_1 =\overline{\mathcal{R}(B^*)}\oplus \mathcal{N}(B)$, we have
$$\mathcal{G}(B,A) =\{(Bx,Ax): x \in \overline{\mathcal{R}(B^*)}\},$$
so if $\mathcal{H}_1=\mathcal{H}_3$  and $\mathcal{R}(B^*)$ is dense in $\mathcal{H}_1$, then $[A/B]$ is densely defined from $\mathcal{H}_1$ to $\mathcal{H}_2$.
\item [(vi)] By Douglas lemma, we have $\mathcal{R}(A^*)\subset \mathcal{R}(B^*) $ if  and only if $\left\|Ax\right\|\leq \mu \left\|Bx\right\|$ for some $\mu >0$ and all $x\in \mathcal{H}_1$. Consequently $[A/B]$ is bounded  from $\mathcal{R}(B)$ to $\mathcal{R}(A)$ if and only if  $\mathcal{R}(A^*)\subset \mathcal{R}(B^*) $.
\item [(vii)] If $H_1 =\mathcal{H}_2=\mathcal{H}_3$, this definition coincides with the one given  by Kaufman, \cite{Kaufman2}.
\item [(viii)] If $H_1 =\mathcal{H}_3$, this definition coincides with the one given by Koliha,  \cite{Koliha}.
\end{itemize}

As remarked  above a right quotient of bounded operators is not necessary bounded neither closed.  The following lemma expresses a  right quotient  in terms of the Moore-Penrose inverse.
\begin{lemma}\label{RQ} Let  $A\in \mathcal{B}(\mathcal{H}_1, \mathcal{H}_2 )$ and $B\in \mathcal{B}(\mathcal{H}_1, \mathcal{H}_3 )$ such that $\mathcal{N}(B) \subseteq \mathcal{N}(A)$. Then $$[A /B]=[AP_{\overline{\mathcal{R}(B^*)}} /B]=AB^{\dag}:  \mathcal{R}(B) \longrightarrow  \mathcal{R}(A).$$
If $\mathcal{R}(B)$ is closed then $[A /B]$ is bounded.
 In particular, if $B$ is invertible then $[A/ B]=AB^{-1}$.
\end{lemma}
\begin{proof} As in the proof of the Lemma \ref{lqd} the Moore-Penrose generalized inverse of $B$ is given by
$$  B^{\dag} =\begin{bmatrix}
 B_1^{-1}& 0\\
 0 & 0
\end{bmatrix}: \begin{bmatrix}
 \mathcal{R}(B)\\
 \mathcal{R}(B)^{\bot}
\end{bmatrix}\longrightarrow  \begin{bmatrix}
 \mathcal{N}(B)^{\bot}\\
 \mathcal{N}(B)
\end{bmatrix} .$$
Since $\mathcal{N}(B) \subseteq \mathcal{N}(A)$, then
$$  A =\begin{bmatrix}
 A_1& 0\\
 A_2 & 0
\end{bmatrix}:\begin{bmatrix}
 \mathcal{N}(B)^{\bot}\\
 \mathcal{N}(B)\end{bmatrix} \longrightarrow \begin{bmatrix}
 \mathcal{R}(B)\\
 \mathcal{R}(B)^{\bot}
\end{bmatrix} 
 .$$

So
$$AB^{\dag}= \begin{bmatrix}
 A_1B_1^{-1}& 0\\
 A_2B_1^{-1}& 0
\end{bmatrix}:   \begin{bmatrix}
\mathcal{R}(B)\\
\mathcal{R}(B)^{\bot}
\end{bmatrix} \longrightarrow \begin{bmatrix}
\mathcal{R}(A)\\
\mathcal{R}(A)^{\bot}
\end{bmatrix}.$$
This implies that
\begin{equation}\label{mpr1}
AB^{\dag}B= \begin{bmatrix}
 A_1& 0\\
A_2 & 0
\end{bmatrix}=A: \begin{bmatrix}
 \mathcal{N}(B)^{\bot}\\
 \mathcal{N}(B)
\end{bmatrix}\longrightarrow   \begin{bmatrix}
 \mathcal{R}(A)\\
\mathcal{R}(A)^{\bot}
\end{bmatrix}.\end{equation}

On the other hand, for all $x\in \mathcal{H}_1 = \mathcal{N}(B)^{\bot}\oplus  \mathcal{N}(B)$, we have
\begin{equation}\label{mpr2}
AB^{\dag}Bx=AP_{\mathcal{N}(B)^{\bot}}x=AP_{\mathcal{N}(B)^{\bot}}x+AP_{\mathcal{N}(B)}x= Ax.
\end{equation}
By \eqref{mpr1} and  \eqref{mpr2}, it follows that  $[A /B]=[AP_{\mathcal{N}(B)^{\bot}} /B]=AB^{\dag}$ with
$$AB^{\dag}: \mathcal{R}(B) \longrightarrow  \mathcal{R}(A).$$
Since $D([A /B])=\mathcal{R}(B)$, it follows that $[A /B]$ is  densely defined in $\overline{\mathcal{R}(B)}$ and is bounded when  $\mathcal{R}(B)$ is closed.
\end{proof}
\begin{remark}For $A\in \mathcal{B}(\mathcal{H})$, we have
	\begin{itemize}
		\item $[A/I]=A$.
		\item If $A$ is injective with closed range, then $[I/A]$  is the left inverse of $A$ and given by $[I/A]=A^{\dag}=(A^*A)^{-1}A^*$.

\end{itemize}
\end{remark}
\begin{corollary}Let $A,  \mathbb{C}^{m\times n}$  and $B\in \mathbb{C}^{p\times n}$ such that $\mathcal{N}(B) \subseteq \mathcal{N}(A)$. Then $[A /B]=AB^{\dag}\in \mathbb{C}^{p\times m}$.
\end{corollary}
If $A\in \mathcal{B}(\mathcal{H}_1, \mathcal{H}_2 )$, $B\in \mathcal{B}(\mathcal{H}_1, \mathcal{H}_3 )$ and $\mathcal{N}(A) = \mathcal{N}(B)$, then we can define $[B / A]$ and $[A / B]$  with
 $$[B / A][A / B] \text{  is the identity operator on }
\mathcal{R}(B) $$ and 
$$[A / B][B / A] \text{  is the identity operator on }
\mathcal{R}(A).$$
 Therefore,  $[A/B]$ is an invertible operator from  $\mathcal{R}(B)$ to $\mathcal{R}(A)$ with
$$ [A / B]^{-1} =[B / A].$$
\begin{corollary}\label{corr1} Let $A\in \mathcal{B}(\mathcal{H}_1, \mathcal{H}_2 )$, $B\in \mathcal{B}(\mathcal{H}_1, \mathcal{H}_3 )$  have the same  kernel. Then $[A/B]$ is an invertible operator from  $\mathcal{R}(B)$ to $\mathcal{R}(A)$ with
$$[A / B]^{-1} =[B / A].$$
If $\mathcal{R}(A)$ is closed then $[A / B]^{-1}$ is bounded.
\end{corollary}

 The graph of a right quotient operator $[A /B]$ (with the condition    $\mathcal{N}(B) \subseteq \mathcal{N}(A)$) is the set
$$\mathcal{G}(B,A) =\{(Bx,Ax): x \in \mathcal{H}_1 \}\subset \mathcal{H}_3  \times \mathcal{H}_2 ,$$
Note that $\mathcal{G}(B,A)$ is the  range of the  operator $T$ defined by:
\begin{equation*}\label{rangeT}
T= \begin{bmatrix}
 B\\
A
\end{bmatrix} : \mathcal{H}_1 \longrightarrow   \mathcal{H}_3  \times \mathcal{H}_2.
\end{equation*}
So its adjoint  is given by
\begin{equation*}
T^*= \begin{bmatrix}
 B^* \quad A^*
\end{bmatrix}: \mathcal{H}_3  \times \mathcal{H}_2 \longrightarrow  \mathcal{H}_1. \end{equation*}
We have 
$\mathcal{R}(T)=\mathcal{G}(B,A)$ and $\mathcal{R}(T^*)=\mathcal{R}(A^*)+\mathcal{R}(B^*)$. Hence
$[A /B]$ is closed if and only if $\mathcal{R}(T)$ is closed in $\mathcal{H}_3  \times \mathcal{H}_2$; equivalently to $\mathcal{R}(T^*)=\mathcal{R}(A^*)+\mathcal{R}(B^*)$ is closed in $\mathcal{H}_1$. This proves the following result.
\begin{lemma}Let  $A \in \mathcal{B}(\mathcal{H}_1, \mathcal{H}_2 )$ and  $ B\in \mathcal{B}(\mathcal{H}_1, \mathcal{H}_3 )$ such that $\mathcal{N}(B) \subseteq \mathcal{N}(A)$. $[A/B]$ is closed if and only if  $\mathcal{R}(A^*)+\mathcal{R}(B^*)$ is closed in $\mathcal{H}_1$.
\end{lemma}
For  $A \in \mathcal{B}(\mathcal{H}_1, \mathcal{H}_2 )$ and  $ B\in \mathcal{B}(\mathcal{H}_1, \mathcal{H}_3 )$, let the following selfadjoint operator
$$R_{A^*,B^*} = (A^*A+B^*B)^{1/2} : \mathcal{H}_1\longrightarrow \mathcal{H}_1;$$
It well known that   $\mathcal{R}(R_{A^*,B^*})=\mathcal{R}(A^*)+\mathcal{R}(B^*)$ \cite[Theorem 2.2]{FilWil71}. Thus 
, if $\mathcal{R}(A^*)+\mathcal{R}(B^*)$ is closed in $\mathcal{H}_1$, then  $\mathcal{N}(R_{A^*,B^*})=\left(\mathcal{R}(A^*)+\mathcal{R}(B^*)\right)^{\perp}= \mathcal{N}(A)\cap \mathcal{N}(B)$. In this case $R_{A^*,B^*}^{\dag}$ exists and bounded, with
$$R_{A^*,B^*}^{\dag}= \begin{bmatrix}
 R & 0\\
 0& 0
\end{bmatrix}:   \begin{bmatrix}
\mathcal{R}(A^*)+\mathcal{R}(B^*)\\
\mathcal{N}(A)\cap \mathcal{N} (B)
\end{bmatrix} \longrightarrow \begin{bmatrix}
\mathcal{R}(A^*)+\mathcal{R}(B^*)\\
\mathcal{N}(A)\cap \mathcal{N} (B)
\end{bmatrix},$$
where $R$ is an isomorphism on  $\mathcal{R}(A^*)+\mathcal{R}(B^*)$. As in \cite{Kaufman1}, we define an operator $J: \mathcal{R}(A^*)+\mathcal{R}(B^*) \longrightarrow \mathcal{H}_3\times \mathcal{H}_2$ by
$$Jx=(BR_{A^*,B^*}^{\dag}x, AR_{A^*,B^*}^{\dag}x).$$
Then for each $x\in \mathcal{R}(A^*)+\mathcal{R}(B^*)$,
\begin{align*}
	\left\|Jx \right\|^2 & = \left\|BR_{A^*,B^*}^{\dag}x\right\|^2+\left\|AR_{A^*,B^*}^{\dag}x\right\|^2\\
	   & = \left\langle x, R_{A^*,B^*}^{\dag}R_{A^*,B^*}^{2}R_{A^*,B^*}^{\dag}x\right\rangle\\
		& = \left\|x \right\|^2
\end{align*}
Hence, if $\mathcal{R}(A^*)+\mathcal{R}(B^*)$ is closed then $J$ is a linear isometry from $\mathcal{R}(A^*)+\mathcal{R}(B^*)$ to $\mathcal{H}_3\times \mathcal{H}_2$. If we assume further, $\mathcal{N}(B) \subseteq \mathcal{N}(A)$.  We have also $ \overline{\mathcal{R}(A^*)}\subset \overline{\mathcal{R}(B^*)}$, this implies that $\mathcal{R}(A^*)+\mathcal{R}(B^*)\subset \overline{\mathcal{R}(B^*)}$. Since $\mathcal{N}(R_{A^*,B^*}) =\mathcal{N}(B)$; then $\mathcal{R}(A^*)+\mathcal{R}(B^*)=\overline{\mathcal{R}(B^*)}$. In this case, since $\mathcal{H}_1 =\overline{\mathcal{R}(B^*)}\oplus \mathcal{N}(B)$, the  graph of the linear operator $[A/B]$is given
$$\mathcal{G}(B,A) =\{(Bx,Ax): x \in \overline{\mathcal{R}(B^*)}=\mathcal{R}(A^*)+\mathcal{R}(B^*)\},$$
and hence the linear isometry $J$ maps the  closed set $ \overline{\mathcal{R}(B^*)}$ onto  $\mathcal{G}(B,A)$. This gives another proof  of the preceding  Lemma.

\section{ Dual properties}
In this section we prove some dual relationships between  the left quotient and the right quotient  of two bounded operators.

Let $A\in \mathcal{B}(\mathcal{H}_1, \mathcal{H}_3 )$ and $B\in \mathcal{B}(\mathcal{H}_2, \mathcal{H}_3 )$. Since $\mathcal{R}(A) \subseteq \mathcal{R}(B)$ then $\mathcal{N}(B^*) \subseteq \mathcal{N}(A^*)$, so we can define $[A^*/B^*]$ the right quotient of $A^*$ by $B^*$ from $\mathcal{R}(B^*)$ to  $\mathcal{R}(A^*)$ such that:
\begin{equation}\label{equ3}
 [A^*/B^*]B^*x=A^*x \qquad \text{ for all } x\in \mathcal{H}_3.
\end{equation}
Since $[B \backslash A]$ is bounded, then its adjoint exists, 
\begin{theorem}\label{dual} Let $A\in \mathcal{B}(\mathcal{H}_1, \mathcal{H}_3 )$ and $B\in \mathcal{B}(\mathcal{H}_2, \mathcal{H}_3 )$  such that $\mathcal{R}(A) \subseteq \mathcal{R}(B)$. We have,  if $\mathcal{R}(B^*)$ is closed, then $$[B \backslash A]^*= [A^*/B^*]=A^*(B^*)^{\dag}.$$
\end{theorem}
\begin{proof}  For all $x,y\in \mathcal{H}_1$,
\begin{align*}
	\left\langle Ax,y\right\rangle & = \left\langle B[B \backslash A] x,y\right\rangle\\
	                               & = \left\langle  x,[B \backslash A]^*B^*y\right\rangle\\
																& = \left\langle  x,A^*y\right\rangle
\end{align*}
Then $$ [B \backslash A]^*B^*y=A^*y\qquad \text{ for all } y\in \mathcal{H}_3.$$
Compare this with \eqref{equ3}, we obtain $[B \backslash A]^*= [A^*/B^*]$. 
\end{proof}
\begin{remark} Without the closure of the range of $B^*$, as seen in the proof of the Douglas lemma, the operator $[B \backslash A]^*$  is the natural extension of  $[A^*/B^*]$ to the closure of the  range of $B^*$; in the following sense.

Since $[A^*/B^*]$ is bounded from $\mathcal{R}(B^*)$ to  $\mathcal{R}(A^*)$; that is $\sup\{ \left\|A^*x\right\|/\left\|B^*x\right\| : x\in \mathcal{H}_3, B^*x\neq 0\}<\infty$, then we can define its natural extension $\widetilde{[A^*/B^*]}$ to $\overline{\mathcal{R}(B^*)}$ by
$$
\widetilde{[A^*/B^*]}x=\left\{ 
	\begin{array}{ll}
		\lim_{n\rightarrow\infty}[A^*/B^*]x_n & \text{ for } x\in \overline{\mathcal{R}(B^*)}, x_n \rightarrow x  \\ 
		0 & \text{ for } x\in \mathcal{R}(B^*)^{\bot}.
	\end{array}
	\right.
$$
\end{remark}

It well known that if $T$ is densely defined from  $\mathcal{H}_1$ to $\mathcal{H}_2$, the adjoint $T^*$  of $T$ exists, is unique and $T^*$ is also densely defined from  $\mathcal{H}_2$ to $\mathcal{H}_1$. Now, let $[A/B]$ be a right quotient with the domain $\mathcal{R}(B)$, witch is  dense in $\overline{\mathcal{R}(B)}$, and let $\mathcal{G}(A,B)$ be its graph. Then the adjoint $G(A,B)^*$ is
naturally defined as the set of elements $(x, y)\in \mathcal{H}_2 \times \mathcal{H}_3  $  such that $\left\langle Az, x\right\rangle =
\left\langle Bz, y\right\rangle$ for all $z\in \mathcal{H}_1$.  We can see that
$$\mathcal{G}(B,A)^*=\{(x,y ):   A^*x=B^*y\}\subset \mathcal{H}_2  \times \mathcal{H}_3 ,$$
and it is a graph again. The corresponding
operator is just the adjoint of $[B/A]$, that is $[B/A]^*$. Furthermore, we have
$$\mathcal{G}(B,A)^*=\mathcal{V}(\mathcal{N}(T^*));$$
with
\begin{equation*}\label{mp1}
T^*= \begin{bmatrix}
 B^* \quad A^*
\end{bmatrix}: \mathcal{H}_3  \times \mathcal{H}_2 \longrightarrow  \mathcal{H}_1 \end{equation*}
and $\mathcal{V}$ is the isomorphisms from  $\mathcal{H}_3  \times \mathcal{H}_2$ to $\mathcal{H}_2 \times \mathcal{H}_3 $ defined by
$$\mathcal{V}(x,y)=(-y,x).$$
\begin{theorem} Let $A \in \mathcal{B}(\mathcal{H}_1, \mathcal{H}_2 )$ and  $ B\in \mathcal{B}(\mathcal{H}_1, \mathcal{H}_3 )$ such that $\mathcal{N}(B) \subseteq \mathcal{N}(A)$. For   a right quotient $[A/B]$, 
\begin{itemize}
	\item its adjoint $[A/B]^*$ exist and closed as an operator from $\overline{R(A)}$ to $\overline{R(B)}$ and
$$[A / B]^*= (B^{\dag})^* A^*;$$
witch uniquely extended  to an operator from $\mathcal{H}_2$ to $\overline{R(B)}$.
\item If $\mathcal{R}(A^*)+\mathcal{R}(B^*)$ is closed in $\mathcal{H}_1$, then $[A/B]^*$ is closed densely defined from $\mathcal{H}_2$ to $\overline{R(B)}$.
\item If $\mathcal{R}(B^*)$ is closed in $\mathcal{H}_1$, then 
$$[A / B]^*= [B^*\backslash A^*]=(B^*)^{\dag}A^*,$$
and
$$[A / B]^{**}= [A / B].$$
\end{itemize}
\end{theorem}
\begin{proof}  Let $[A/B]$  a right quotient with  domain $\mathcal{R}(B)$, witch is  dense in $\overline{\mathcal{R}(B)}$, then its adjoint $[A/B]^*$ exist and closed as an operator from $\overline{R(A)}$ to $\overline{R(B)}$ with
$$ \mathcal{D}([A/B]^*)=\{ x\in \overline{R(A)} : \exists y \in \overline{R(B)} , \left\langle [A/B]Az, x\right\rangle =
\left\langle Bz, y\right\rangle ; \forall  z\in \mathcal{H}_1\}$$
and
$$[A/B]^*x=y. $$
 By lemma \ref{RQ} $[A/B]=AB^{\dag}$, this implies that $[A/B]^*=(B^{\dag})^* A^*.$
Further, for all $x\in \mathcal{D}([A/B]^*)$, 
$$B^*[A/B]^*x= B^*(B^{\dag})^* A^*x=(B^{\dag}B)^* A^*x= P_{\overline{\mathcal{R}(B^*)}}A^*x =A^*x;$$ 
this holds because $\mathcal{R}(A^*)\subset \overline{\mathcal{R}(B^*)}$. The last equality is also valid if $x\in  \mathcal{N}(A^*)$, so is holds for all $x\in \mathcal{D}([A/B]^*)\oplus \mathcal{N}(A^*)\subseteq \mathcal{H}_2$; this implies that we can extended $[A/B]^*$ to a closed operator from $\mathcal{H}_2$ to $\overline{R(B)}$. Further; if $\mathcal{R}(A^*)+\mathcal{R}(B^*)$ is closed in $\mathcal{H}_1$, then 
$[A/B]$ is closed densely defined from $\overline{R(B)}$ to $\overline{R(A)}$. So, $[A/B]^*$ is closed densely defined from $\mathcal{H}_2$ to $\overline{R(B)}$, with $\mathcal{D}([A/B]^*)$ is dense in $\overline{\mathcal{R}(A)}$.

Now if   $\mathcal{R}(B^*)$ is closed in $\mathcal{H}_1$, in that case  $\mathcal{R}(A^*)\subset \mathcal{R}(B^*)$, so $[A / B]$ is bounded and  its adjoint $[A/B]^*$ exist and bounded. Again by lemma \ref{RQ}, 
$$[A/B]^*=(B^{\dag})^* A^*= (B^*)^{\dag} A^*=[B^*\backslash A^*],$$
and 
$$ B^*[B^*\backslash A^*]x=B^* (B^*)^{\dag}A^*x= P_{\mathcal{R}(B^*)}A^*x =A^*x \quad \text{ for all } x\in \mathcal{H}_2.$$ 
Then $$ B^*y=A^*x\qquad \text{ for all } x\in \mathcal{H}_2 \quad \text{ and } \quad y=[B^*\backslash A^*]x \in \mathcal{H}_3.$$
Hence $[A / B]^*= [B^*\backslash A^*]=(B^*)^{\dag}A^*$ and $[A / B]^{**}= [A / B].$
\end{proof}
\begin{remark} Let us point out that Izumino, in \cite{izumino89}, construct  a right quotient operator to be an adjoint of  $[A/B]$, so our objective is completely different. Here the adjoint  of  $[A/B]$ is always exists (in particular, as a left quotient) but its adjoint in the same nature (as a right quotient) may be not exist, for more details see \cite{izumino89}.
\end{remark}

\section{Sums and products}
Let $A$, $B$, $C$, $D$ $\in \mathcal{B}(\mathcal{H}_1, \mathcal{H}_2 )$. Let $[B\backslash A]$  and $[D\backslash C]$ be two left quotient operators  (with the respective range conditions). Then we prove that the sum of them is also a left quotient operators. To define the denominator of the sum we need the following range operator, called parallel sum of $B^*B$ and $D^*D$, defined  by
$$ B^*B : D^*D = B^*B(B^*B+D^*D)^{\dag}D^*D : \mathcal{H}_1 \longrightarrow \mathcal{H}_1 $$
with $\mathcal{R}(B^*)+ \mathcal{R}(D^*)$ is closed in $\mathcal{H}_1$, see \cite[p.p. 277]{FilWil71}.
If we denotes $$S_{B^*,D^*}=(B^*B : D^*D)^{1/2}$$ then we have, by \cite[Theorem  4.2]{FilWil71}:
\begin{enumerate}
	\item $B^*B : D^*D=  D^*D: B^*B$ is positive and bounded.
	\item $\mathcal{R}(S_{B^*,D^*})=\mathcal{R}(B^*)\cap \mathcal{R}(D^*)$.
\end{enumerate}

 This implies that  the left operators quotient
$$  B_1^*=[B^* \backslash S_{B^*,D^*}] \qquad  \text{ and } \qquad  D_1^*=[D^* \backslash S_{B^*,D^*}] $$
and the right operators quotient
$$ B_1=[ S_{B^*,D^*} / B] \qquad  \text{ and } \qquad D_1=[ S_{B^*,D^*} / D]$$
are always makes sens. 
Now, for the left quotient form of the sum we have
\begin{theorem} Let $A$, $B$, $C$, $D$ $\in \mathcal{B}(\mathcal{H}_1, \mathcal{H}_2 )$. Let $[B\backslash A]$  and $[D\backslash C]$ be two left quotient operators  (with the respective range conditions). We have
$$[B\backslash A]+[D\backslash C]=[S_{B^*,D^*}\backslash (B_1A+D_1C)] : \mathcal{H}_1 \longrightarrow \mathcal{N}(B)^{\bot} + \mathcal{N}(D)^{\bot}.$$
\end{theorem}
\begin{proof} If  $\mathcal{R}(A)\subseteq \mathcal{R}(B)$ and $\mathcal{R}(C)\subseteq \mathcal{R}(D)$, then the operators $B_1A+D_1C$ is well-defined from  $\mathcal{H}_1$ to $\mathcal{R}(S_{B^*,D^*})$. Hence $\mathcal{R}(B_1A+D_1C)\subset \mathcal{R}(S_{B^*,D^*})$ and $[S_{B^*,D^*}\backslash B_1C+D_1A]$ is also well-defined from  $\mathcal{H}_1$ to $ \overline{\mathcal{R}(S_{A^*,B^*})}$. Further,
$$[S_{B^*,D^*}\backslash ( B_1A+D_1C)]=S_{B^*,D^*}^{\dag}(B_1A+ D_1C).$$
On the other hand, $B_1=S_{B^*,D^*}P_{\overline{\mathcal{R}(B^*)}}B^{\dag}$ and $D_1=S_{B^*,D^*}P_{\overline{\mathcal{R}(D^*)}}D^{\dag}$, so
\begin{align*}
	S_{B^*,D^*}^{\dag}(B_1A+ D_1C) & = S_{B^*,D^*}^{\dag}(S_{B^*,D^*}P_{\overline{\mathcal{R}(B^*)}}B^{\dag}A+S_{B^*,D^*}P_{\overline{\mathcal{R}(D^*)}}D^{\dag}C) \\
	                               & =  P_{\overline{\mathcal{R}(S_{B^*,D^*})}}(P_{\overline{\mathcal{R}(B^*)}}B^{\dag}A+P_{\overline{\mathcal{R}(D^*)}}D^{\dag}C)\\
																& =P_{\overline{\mathcal{R}(S_{B^*,D^*})}}([B\backslash A]+[D\backslash C]).
\end{align*}
On the other hand; it is clear that the sum $[B\backslash A]+[D\backslash C]$ is an operator defined from $\mathcal{H}_1$ to $\mathcal{N}(B)^{\bot} + \mathcal{N}(D)^{\bot}= \overline{\mathcal{R}(S_{B^*,D^*})}$, which implies that
$$ [B\backslash A]+[D\backslash C]=[S_{B^*,D^*}\backslash B_1A+D_1C].$$
\end{proof}
\begin{remark} In the case of the scalars case, that is $a, b, c, d \in \mathbb{C}$, this Theorem  is nothing,
$$\frac{a}{b} +\frac{c}{d}=\frac{b_1 a+d_1 c}{\left|b\right| \left(\sqrt{\left|b\right|^2 +\left|d\right|^2}\right)^{-1}\left|d\right|}, $$with $$ b_1=\frac{\sqrt{\left|b\right|^2 +\left|d\right|^2}}{b\left|d\right|\left|b\right|}, \quad d_1=\frac{\sqrt{\left|b\right|^2 +\left|d\right|^2}}{d\left|d\right|\left|b\right|} \quad b\neq 0, b\neq 0.$$
\end{remark}
For two  left quotient operators with the same denominator, we have,

\begin{theorem}\label{thmsomme} Let $A$, $B$ $\in \mathcal{B}(\mathcal{H}_1, \mathcal{H}_3 )$ and $D\in \mathcal{B}(\mathcal{H}_2, \mathcal{H}_3 )$ such that $\mathcal{R}(A) \subseteq \mathcal{R}(D)$ and $\mathcal{R}(B)\subseteq \mathcal{R}(D)$. Then
$$ [D \backslash (A+B)]=[D \backslash A]+[D \backslash B].$$
\end{theorem}
\begin{proof}By assumptions, we have $A=D[D \backslash A]$ and $B=D[D \backslash B]$, so 
$$ A+B=D([D \backslash A]+[D \backslash B]).$$
On the other hand; $\mathcal{R}(A+B) \subseteq \mathcal{R}(A)+\mathcal{R}(B) \subseteq \mathcal{R}(D)$, it follows that $[D \backslash A+B]$ is well defined and
$$ A+B=D[D \backslash (A+ B)].$$
By what we assert the desired equality.
\end{proof}
We next show that the product of two left quotient operators is a gain a left operator.
\begin{theorem}\label{produitleft} Let $A$, $B$ $\in \mathcal{B}(\mathcal{H}_1, \mathcal{H}_2 )$ and $C$, $D$ $\in \mathcal{B}(\mathcal{H}_1, \mathcal{H}_3 )$. If $[B\backslash A]$  and $[D\backslash C]$ are two left quotient operators  (with the respective range conditions), then
$$[B\backslash A][D\backslash C]=[NB\backslash MC],$$
where $M$ $\in \mathcal{B}(\mathcal{H}_3, \mathcal{H}_1 )$ and $N$ $\in \mathcal{B}(\mathcal{H}_2, \mathcal{H}_1 )$ satisfying the conditions  $MD=NA$ and $\mathcal{N}(N)^{\bot}=\mathcal{R}(B)$.
\end{theorem}
To prove this theorem, we need the following lemma
\begin{lemma}\label{productPM} If $T$ $\in \mathcal{B}(\mathcal{H}_1, \mathcal{H}_2 )$ and $S$ $\in \mathcal{B}(\mathcal{H}_2, \mathcal{H}_3 )$, with $\mathcal{R}(T)=\mathcal{N}(S)^{\bot}$; then $(ST)^{\dag}=T^{\dag}S^{\dag}$.
\end{lemma}
\begin{proof} In this case we have $\mathcal{R}(T)$ is closed and $ ST\in \mathcal{B}(\mathcal{H}_1, \mathcal{H}_3 )$ with
  $$(ST)^{\dag}= \begin{bmatrix}
 R_1^{-1}& 0\\
 0& 0
\end{bmatrix}:   \begin{bmatrix}
\mathcal{R}(ST)\\
\mathcal{R}(ST)^{\perp}
\end{bmatrix} \longrightarrow \begin{bmatrix}
\mathcal{N}(ST)^{\perp}\\
\mathcal{N}(ST)
\end{bmatrix}.$$
Now, if 
$$T^{\dag}= \begin{bmatrix}
 T_1^{-1}& 0\\
 0& 0
\end{bmatrix}:   \begin{bmatrix}
\mathcal{R}(T)\\
\mathcal{R}(T)^{\perp}
\end{bmatrix} \longrightarrow \begin{bmatrix}
\mathcal{N}(T)^{\perp}\\
\mathcal{N}(T)
\end{bmatrix}$$
and
$$S^{\dag}= \begin{bmatrix}
 S_1^{-1}& 0\\
 0& 0
\end{bmatrix}:   \begin{bmatrix}
\mathcal{R}(S)\\
\mathcal{R}(S)^{\perp}
\end{bmatrix} \longrightarrow \begin{bmatrix}
\mathcal{N}(S)^{\perp}\\
\mathcal{N}(S)
\end{bmatrix}$$
then
$$T^{\dag}S^{\dag}= \begin{bmatrix}
 T_1^{-1}S_1^{-1}& 0\\
 0& 0
\end{bmatrix}:   \begin{bmatrix}
\mathcal{R}(S)\\
\mathcal{R}(S)^{\perp}
\end{bmatrix} \longrightarrow \begin{bmatrix}
\mathcal{N}(T)^{\perp}\\
\mathcal{N}(T)
\end{bmatrix};$$
taking in account that $\mathcal{N}(S)^{\bot}=\mathcal{R}(T)$. We have $\mathcal{R}(ST)\subset \mathcal{R}(S)$. Now; let $y\in \mathcal{R}(S)=\mathcal{R}(SS^*)$, so $y=SS^*x$ for $x \in \mathcal{H}_3$. By assumption $\mathcal{R}(S^*)=\mathcal{N}(S)^{\bot}=\mathcal{N}(S^*S)^{\bot}=\mathcal{R}(S^*S)=\mathcal{R}(T)=\mathcal{R}(TT^*)$. Thus $S^*x\in \mathcal{R}(TT^*)$ and $S^*x=TT^*z$, for $z\in  \mathcal{H}_3$. So $y=STT^*z$ witch implies that $y\in \mathcal{R}(ST)$. Consequently, $\mathcal{R}(ST)= \mathcal{R}(S)$. This proves that $(ST)^{\dag}$ and $T^{\dag}S^{\dag}$ have the same domain and $\mathcal{N}(T^{\dag}S^{\dag})^{\bot})=\mathcal{N}((ST)^{\dag})=\mathcal{R}(S)^{\bot}$. Let us remark that $TT^{\dag}=P_{\mathcal{R}(T)}=P_{\mathcal{N}(S)^{\bot}}= SS^{\dag}$.  For all $y \in \mathcal{H}_3$, 
\begin{align*}
	T^{\dag}S^{\dag}STT^{\dag}S^{\dag}y&=T^{\dag}S^{\dag}SS^{\dag}SS^{\dag}y\\
	&=T^{\dag}P_{\overline{\mathcal{R}(S^{\dag})}}S^{\dag}y\\
	&=T^{\dag}S^{\dag}y. 
\end{align*}
 Hence  $T^{\dag}S^{\dag}STx=P_{\overline{\mathcal{R}(T^{\dag}S^{\dag})}}x$, for all $y \in \mathcal{H}_1$.

Similarly, for all $x \in \mathcal{H}_1$, 
\begin{align*}
ST	T^{\dag}S^{\dag}STx&=SS^{\dag}SS^{\dag}STx\\
	&=SP_{\overline{\mathcal{R}(S)}}Tx\\
	&=SP_{\overline{\mathcal{R}(ST)}}Tx\\
	&=STx. 
\end{align*}
 Hence  $STT^{\dag}S^{\dag}z=P_{\overline{\mathcal{R}(ST)}}z$, for all $z \in \mathcal{D}(T^{\dag}S^{\dag})=\mathcal{R}(S)\oplus \mathcal{R}(S)^{\perp}$.
Hence, by the uniqueness of Moore-Penrose inverse, $(ST)^{\dag}=T^{\dag}S^{\dag}$.
\end{proof}
\begin{proof}[Proof of Theorem \ref{produitleft}] We have $[B\backslash A][D\backslash C]=B^{\dag}AD^{\dag}C$ is an operator well defined from $ \mathcal{H}_1$ to $\mathcal{N}(B)^{\bot}$ and the operator $[NB\backslash MC]$ is defined from $ \mathcal{H}_1$ to $\mathcal{N}(NB)^{\bot}\subset \mathcal{N}(B)^{\bot}$. Now, by assumption, $MD=NA$, so $N^{\dag}MDD^{\dag}=N^{\dag}NAD^{\dag}$; Thus  $N^{\dag}MP_{\overline{\mathcal{R}(D)}}=P_{\overline{\mathcal{N}(N)^{\perp}}}AD^{\dag}=P_{\overline{\mathcal{R}(B)}}AD^{\dag}=AD^{\dag}$ ($\mathcal{R}(A)\subset \mathcal{R}(B)$). This implies that $B^{\dag}AD^{\dag}C=B^{\dag}N^{\dag}MP_{\overline{\mathcal{R}(D)}}C=B^{\dag}N^{\dag}MC$. Now by Lemma \ref{productPM}, $B^{\dag}N^{\dag}=(NB)^{\dag}$; then 
$$[B\backslash A][D\backslash C]=B^{\dag}AD^{\dag}C=(NB)^{\dag}MC=[NB\backslash MC].$$
\end{proof}
Now using Theorem \ref{dual}, we get the following results for the sum and product for the right quotient operators.

\begin{theorem}\label{thmsommeR} Let $A$, $B$, $C$, $D$ $\in \mathcal{B}(\mathcal{H}_1, \mathcal{H}_2 )$. If $[A/ B]$  and $[C/ D]$ are two right quotient operators  (with the respective kernel conditions) sucht that $\mathcal{R}(B)+ \mathcal{R}(D)$ is closed in $\mathcal{H}_2$, then
$$[A/ B]+[C/ D]=[ AB_1+CD_1/S_{B,D}]: \mathcal{R}(B)\cap \mathcal{R}(D)  \longrightarrow \mathcal{R}(A)\cap \mathcal{R}(C),$$
with 
$$ B_1=[ B \backslash  S_{B,D} ] \qquad  \text{ and } \qquad D_1=[ D\backslash  S_{B,D} ].$$
\end{theorem}

\begin{theorem} Let $A$, $B$ $\in \mathcal{B}(\mathcal{H}_1, \mathcal{H}_2 )$ and $D\in \mathcal{B}(\mathcal{H}_1, \mathcal{H}_3 )$ such that $\mathcal{N}(D) \subseteq \mathcal{N}(A)$ and $\mathcal{N}(D)\subseteq \mathcal{N}(B)$. Then
$$ [(A+B)/D ]=[A/ D  ]+[ B/D].$$
\end{theorem}

\begin{theorem}\label{thmproductR} Let $A$, $B$ $\in \mathcal{B}(\mathcal{H}_1, \mathcal{H}_3 )$ and $C$, $D$ $\in \mathcal{B}(\mathcal{H}_2, \mathcal{H}_3 )$. If $[ A/B]$  and $[C/D]$  are two right quotient operators  (with the respective kernel conditions);  then
$$[A/ B][C/ D]=[AP_{\mathcal{N}(B)^{\bot}} M/DN],$$
where $M$ $\in \mathcal{B}(\mathcal{H}_3, \mathcal{H}_1 )$ and $N$ $\in \mathcal{B}(\mathcal{H}_3, \mathcal{H}_2 )$ satisfying the conditions   $BM=CN$  and $\mathcal{N}(D)^{\bot}=\mathcal{R}(N)$. 
\end{theorem}
\begin{proof} Since the domain of $[A/ B][C/ D]$ is $\mathcal{R}(D)$ and those of $[AM/DN]$ is $\mathcal{R}(DN)$, as seen in Lemma \ref{productPM}, the assumption $\mathcal{N}(D)^{\bot}=\mathcal{R}(N)$ guaranties that this two operators have the same domain. Now, by assumption, $BM=CN$, so $B^{\dag}BMN^{\dag} = B^{\dag}CNN^{\dag}$; Thus  $P_{\mathcal{N}(B)^{\perp}}MN^{\dag}=B^{\dag}C P_{\mathcal{N}(D)^{\perp}}$. This implies that for $y\in \mathcal{R}(D)$, we have
\begin{align*}
	[A/ B][C/ D]y& =AB^{\dag}CD^{\dag}y=AB^{\dag}CP_{\mathcal{N}(D)^{\perp}}D^{\dag}y\\
	&=AP_{\mathcal{N}(B)^{\perp}}MN^{\dag}D^{\dag}y\\
	&=AMN^{\dag}D^{\dag}y\\
	&=AM(DN)^{\dag}y\\
	&= [AM/DN]y.
\end{align*}
Here we use the fact that $AP_{\mathcal{N}(B)^{\perp}}=A$ because $\mathcal{N}(B)\subset \mathcal{N}(A) $ and $(DN)^{\dag}=N^{\dag}D^{\dag}$ by Lemma \ref{productPM}.
\end{proof}
\begin{remark} Theorems \ref{thmsommeR} and \ref{thmproductR} are a generalization of 
\cite[ Theorem 3.1]{izumino89} and  \cite[ Theorem 3.2]{izumino89}, respectively.
\end{remark}
We finally turn to the Lemma \ref{productPM}, to give some results like a simplification in a quotient as in the case of scalars.
\begin{corollary} Let $A \in \mathcal{B}(\mathcal{H}_1, \mathcal{H}_3 )$,  $ B\in \mathcal{B}(\mathcal{H}_2, \mathcal{H}_3 )$ and $[B\backslash A]$  is a left quotient operator, we have
$$ [MB \backslash MA]=[ B \backslash A]$$
for all $M\in \mathcal{B}(\mathcal{H}_3, \mathcal{H}_4 )$ such that $\mathcal{N}(M)^{\bot}=\mathcal{R}(B)$. In particular, if $M\in \mathcal{B}(\mathcal{H}_3, \mathcal{H}_2 )$ and $M=B^*$, we have
$$ [B^*B \backslash B^*A]=[ B \backslash A].$$
\end{corollary}
\begin{proof}  We have $\mathcal{R}(AM)\subset \mathcal{R}(A)\subset \mathcal{R}(B)=\mathcal{R}(BM)$; then $[MB \backslash MA]$ is well defined. Further,
\begin{align*}
	[MB \backslash MA]&=(MB)^{\dag}MA=B^{\dag}M^{\dag}MA\\
	&=B^{\dag}P_{\mathcal{N}(M)^{\perp}}A
	=B^{\dag}P_{\mathcal{R}(B)}A=B^{\dag}A=[ B \backslash A].
\end{align*}
\end{proof}
\begin{corollary} Let $A \in \mathcal{B}(\mathcal{H}_1, \mathcal{H}_2 )$,  $ B\in \mathcal{B}(\mathcal{H}_1, \mathcal{H}_3 )$ and $[ A/B]$  is a right quotient operator, we have
$$ [ AM/BM]=[ A/B]$$
for all $M\in \mathcal{B}(\mathcal{H}_4, \mathcal{H}_1 )$ such that $\mathcal{N}(B)^{\bot}=\mathcal{R}(M)$. In particular, if $M\in \mathcal{B}(\mathcal{H}_3, \mathcal{H}_1 )$ and $M=B^*$, we have
$$ [ AB^*/BB^*]=[ A/B].$$
\end{corollary}
\begin{proof}  Since $\mathcal{N}(B)^{\bot}=\mathcal{R}(M)$ and $\mathcal{N}(B)\subset \mathcal{N}(A)$; then $ \mathcal{N}(BM)\subset\mathcal{N}(AM)$; then $[ AM/BM]$ is well defined. Further,
\begin{align*}
	[ AM/BM] & =AM(BM)^{\dag}=AMM^{\dag}B^{\dag}\\
	& =AP_{\mathcal{R}(M)}B^{\dag}=AP_{\mathcal{N}(B)^{\bot}}B^{\dag}=AB^{\dag}=[ A/B].
\end{align*}
\end{proof}
\begin{corollary} Let $A \in \mathcal{B}(\mathcal{H}_1, \mathcal{H}_2 )$ with closed range, we have
$$ A=[ A^*\backslash A^*A] \quad \text{ and } \quad A^{\dag}=[ A^*A\backslash A^*].$$
\end{corollary}


\end{document}